\documentclass[twoside,12pt]{article}
\usepackage[]{amsmath,amsthm,amssymb}
\usepackage[latin1]{inputenc}

\begin{document}

\title{{\Large\bf  Discrete Lebedev-Skalskaya transforms}}

\author{Semyon  YAKUBOVICH}
\maketitle

\markboth{\rm \centerline{ Semyon  YAKUBOVICH}}{}
\markright{\rm \centerline{LEBEDEV-SKALSKAYA TRANSFORMS}}

\begin{abstract} {\noindent Discrete analogs of the Lebedev-Skalskaya transforms are introduced and investigated. It involves series and integrals with respect to the kernels ${\rm Re} K_{\alpha+in}(x),  {\rm Im} K_{\alpha+in}(x),\\ x >0, n \in \mathbb{N}, |\alpha | < 1,\  i $ is the imaginary unit and $K_\nu(z)$ is the modified Bessel function. The corresponding inversion formulas for suitable  functions and sequences in terms of these series and integrals are established when $\alpha = \pm 1/2$.  The case $\alpha=0$  reduces to the Kontorovich-Lebedev transform.}

\end{abstract}
\vspace{4mm}

{\bf Keywords}: {\it Lebedev-Skalskaya  transforms,   Kontorovich-Lebedev transform, modified Bessel function, Fourier series}

{\bf AMS subject classification}:  45A05,  44A15,  42A16, 33C10

\vspace{4mm}

\section {Introduction and preliminary results}

This investigation deals  with discrete analogs of the following reciprocal pairs of the Lebedev-Skalskaya  transformations

$$F(\tau)= \int_0^\infty  {\rm Re} K_{1/2+i\tau}(x) f(x) dx,\quad \tau \in \mathbb{R},\eqno(1.1)$$

$$f(x)= {4\over \pi^2} \int_{0}^\infty \cosh(\pi\tau)  {\rm Re} K_{ 1/2+i\tau}(x) F(\tau) d\tau,\quad x > 0,\eqno(1.2)$$

$$F(\tau)= \int_0^\infty  {\rm Im} K_{ 1/2+i\tau}(x) f(x) dx, \quad \tau \in \mathbb{R},\eqno(1.3)$$

$$f(x)= {4\over \pi^2} \int_{0}^\infty  \cosh(\pi\tau)  {\rm Im} K_{1/2+i\tau}(x) F(\tau) d\tau,\quad x > 0.\eqno(1.4)$$
Here ${\rm Re},\ {\rm Im}$ are real and imaginary parts, respectively,  of a complex-valued function and $K_\nu(z)$ is the modified Bessel function [4] which is represented by the integral

$$ K_\nu(z) = \int_0^\infty e^{-x\cosh(u)} \cosh(\nu u) du,\  {\rm Re} z >0,\  \nu \in \mathbb{C}.\eqno(1.5)$$
The transformations (1.1)-(1.4) were considered for the first time in [2] (cf. [4]) in connection with applications to the related problems in mathematical physics. Formula  (1.5)  immediately suggests  integral representations for the functions ${\rm Re} K_{\alpha+i\tau}(x),  {\rm Im} K_{\alpha+i\tau}(x),\ \alpha \in \mathbb{R}$

$$ {\rm Re} K_{\alpha+i\tau}(x) =  \int_0^\infty e^{-x\cosh(u)} \cosh(\alpha u) \cos(\tau u) du,\eqno(1.6)$$

$$ {\rm Im} K_{\alpha+i\tau}(x) =  \int_0^\infty e^{-x\cosh(u)} \sinh(\alpha u) \sin(\tau u) du.\eqno(1.7)$$
These functions satisfy the following inequality (see [4], p. 172)

$$\left| \left\{\begin{aligned} {\rm Re} \\ {\rm Im} \end{aligned}\right\}  K_{\alpha + i\tau} (x) \right| \le e^{-\delta |\tau|} K_\alpha (x\cos(\delta)),\quad \delta \in \left[0, {\pi\over 2}\right),\ x >0, \ \tau \in \mathbb{R}.\eqno(1.8)$$
The case $\alpha =0$ corresponds to the real-valued function $K_{i\tau}(x)$ which is the kernel of the Kontorovich-Lebedev transform [4] whose discrete analogs were considered recently by the author in [5].  Our goal here is to investigate and prove inversion formulas for $\alpha=\pm 1/2$ in suitable classes of sequences and functions for the following discrete Lebedev-Skalskaya transformations

$$f(x)= \sum_{n=0}^\infty a_n \ {\rm Re} K_{\alpha+in}(x),\quad x >0,\eqno(1.9)$$

$$f(x)= \sum_{n=0}^\infty a_n \ {\rm Im} K_{\alpha+in}(x),\quad x >0,\eqno(1.10)$$

$$a_n= {4\over \pi^2}  \int_0^\infty  {\rm Re} K_{\alpha+in }(x) f(x) dx,\quad n \in \mathbb{N}_0,\eqno(1.11)$$

$$a_n=  {4\over \pi^2}  \int_0^\infty  {\rm Im} K_{\alpha+in }(x) f(x) dx,\quad  n \in \mathbb{N}_0.\eqno(1.12)$$
To do this we will follow in the sequel  the technique developed in [5] for the discrete Kontorovich-Lebedev transform, employing  the following integral (see [3], Vol. II, Entry 2.16.6.1)

$$ \int_0^\infty e^{-x\cosh(u)}  K_{\alpha+i\tau}(x) dx =  {\pi \sinh ( (\alpha+i\tau)  u) \over \sinh (u) \sin( (\alpha+i\tau)\pi)}, \quad u, \tau \in \mathbb{R}, \  |\alpha| < 1.\eqno(1.13)$$
Hence we find $(n \in \mathbb{N})$

$$ \int_0^\infty e^{-x\cosh(u)}  {\rm Re} K_{\alpha+in}(x) dx  $$

$$=   {\pi   [ \sinh(\alpha u)  \sin(\alpha\pi)  \cosh(\pi n) \cos(nu)    + \cosh(\alpha u) \cos(\alpha\pi) \sinh(\pi n) \sin(nu)   ] \over \sinh(u)  [\sinh^2(\pi n) + \sin^2(\pi \alpha ) ]},\eqno(1.14) $$

$$ \int_0^\infty e^{-x\cosh(u)}  {\rm Im} K_{\alpha+in}(x) dx  $$

$$=   {\pi   [ \cosh(\alpha u)  \sin(\alpha\pi)  \cosh(\pi n) \sin(nu)    -  \sinh(\alpha u) \cos(\alpha\pi) \sinh(\pi n) \cos(nu)   ] \over \sinh(u)  [\sinh^2(\pi n) + \sin^2(\pi \alpha) ]}.\eqno(1.15) $$
Finally in this section we introduce,  following [1],   the incomplete modified Bessel function $J(z,\nu,w)$

$$J(z,\nu,w) = \int_{0}^w  e^{-z\cosh (u) } \cosh(\nu u) du.\eqno(1.16)$$ 
Therefore, in particular,  this formula will immediately imply 

$$ {\rm Re} J(x, \alpha+i\tau, \pi) =  \int_0^\pi e^{-x\cosh(u)} \cosh(\alpha u) \cos(\tau u) du,\  x > 0,\eqno(1.17)$$

$$ {\rm Im} J(x, \alpha+i\tau, \pi)  =  \int_0^\pi  e^{-x\cosh(u)} \sinh(\alpha u) \sin(\tau u) du,\ x > 0.\eqno(1.18)$$

\section{Inversion and representation theorems}

We begin with

{\bf Theorem 1}. {\it Let $\alpha \in \mathbb{R},\ |\alpha| < 1.$ Let the sequence $\{a_n\}_{n\ge 0}$ be such that  

$$ \sum_{n=0}^\infty \left| a_n \right| e^{-\delta n}    < \infty,\quad  \delta \in \left[0, \  {\pi\over 2} \right).\eqno(2.1)$$
Then compositions of Laplace transform and transformations $(1.9), (1.10)$ have the following representations, respectively,

$$ \int_0^\infty e^{-x\cosh(u)} f(x) dx = { \pi \over \sinh(u)}  \left[  \sum_{n=0}^\infty {a_n \over \sinh^2(\pi n) + \sin^2(\pi \alpha)} \right.$$

$$\Bigg.\left[ \sinh(\alpha u)  \sin(\alpha\pi)  \cosh(\pi n) \cos(nu)  +  \cosh(\alpha u) \cos(\alpha\pi) \sinh(\pi n) \sin(nu) \right] \Bigg],\eqno(2.2)$$
where $ u \in \mathbb{R}$,

$$ \int_0^\infty e^{-x\cosh(u)} f(x) dx = { \pi \over \sinh(u)}  \left[  \sum_{n=1}^\infty {a_n \over \sinh^2(\pi n) + \sin^2(\pi \alpha)} \right.$$

$$\Bigg.\left[ \cosh(\alpha u)  \sin(\alpha\pi)  \cosh(\pi n) \sin(nu)    -  \sinh(\alpha u) \cos(\alpha\pi) \sinh(\pi n) \cos(nu)\right] \Bigg],\eqno(2.3)$$
where $u \in \mathbb{R}, a_0=0$, series and integrals in  $(2.2), (2.3)$ converge absolutely.}

\begin{proof} The proof is straightforward from (1.8),  (2.1)   and asymptotic behavior of the modified Bessel function [4], which guarantee the estimate

$$\int_0^\infty  e^{-x\cosh(u)} \left|  f(x)\right| dx \le  \sum_{n=0}^\infty |a_n|  e^{-\delta n}   \int_0^\infty  e^{-x}  K_{\alpha}(x \cos(\delta)) dx < \infty\eqno(2.4)$$
when $ |\alpha| < 1,   \delta \in \left[0, \  \pi/ 2 \right)$.   Therefore, plugging formulas (1.9), (1.10) for $f$ in the left-hand side of equalities (2.2), (2.3), respectively, we  interchange the order of integration and summation by Fubini's  theorem and then appeal to (1.14), (1.15) to get the desired equalities (2.2), (2.3). 
\end{proof}

The inversion formulas for discrete transformations (1.9), (1.10) when $\alpha = \pm 1/2$ are established by the following theorem.

{\bf Theorem 2}. {\it Let $\alpha= \pm 1/2$ and the sequence $\{a_n\}_{n\ge 0}$ satisfy conditions of Theorem $1$. Then transformations $(1.9), (1.10)$ can be inverted by the formulas, respectively,

$$a_n= {4\over \pi^2} \cosh(\pi n) \int_0^\infty  {\rm Re} J(x, \pm 1/2+in, \pi) f(x) dx,\eqno(2.5)$$

$$a_n= {4\over \pi^2} \cosh(\pi n) \int_0^\infty  {\rm Im} J(x, \pm 1/2+in, \pi) f(x) dx,\eqno(2.6)$$
where functions $ {\rm Re} J(x, \alpha+i\tau, \pi),  {\rm Im} J(x,  \alpha+i\tau, \pi)$ are defined by $(1.17), (1.18)$ and integrals $(2.5), (2.6)$ converge absolutely.}

\begin{proof} Since, evidently,  ${\rm Re} K_{-1/2+in}(x) = {\rm Re} K_{1/2+in}(x),\   {\rm Im} K_{-1/2+in}(x) = - {\rm Im} K_{1/2+in}(x)$ via (1.5), it is sufficient to prove the theorem for $\alpha =1/2$, which corresponds to the case of the Lebedev-Skalskaya transforms (1.1), (1.3).  To do this we write (2.2), (2.3) accordingly,

$$ \int_0^\infty e^{-x\cosh(u)} f(x) dx = { \pi   \over 2\cosh(u/2)}  \sum_{m=0}^\infty {a_m  \cos(mu) \over \cosh(\pi m)},\eqno(2.7)$$ 

$$ \int_0^\infty e^{-x\cosh(u)} f(x) dx = { \pi  \over 2\sinh(u/2)}  \sum_{m=0}^\infty {a_m  \sin(mu) \over \cosh(\pi m)},\eqno(2.8)$$ 
taking into account $f$ by formulas (1.9), (1.10), correspondingly.  Hence for the right-hand side of (2.5)  we derive the equalities via  (1.17), (2.7) 

$${4\over \pi^2} \cosh(\pi n) \int_0^\infty  {\rm Re} J(x,  1/2+in, \pi) f(x) dx $$

$$=  {4\over \pi^2} \cosh(\pi n) \int_0^\infty  \int_0^\pi e^{-x\cosh(u)} \cosh( u/2) \cos(n u) du f(x) dx $$

$$=  {2\over \pi} \cosh(\pi n)  \int_0^\pi  \cos(n u)  \sum_{m=0}^\infty {a_m  \cos(mu) \over \cosh(\pi m)} du = a_n,$$
where the interchange of the order of integration and summation is allowed by Fubini's theorem owing to the estimate

$$\int_0^\infty \left| {\rm Re} J(x,  1/2+in, \pi) f(x)\right| dx \le   \sum_{m=0}^\infty |a_m|  e^{-\delta m}  \int_0^\pi \cosh(u/2) du$$

$$\times   \int_0^\infty  e^{-x}  K_{1/2}(x \cos(\delta)) dx = {\pi \sqrt{2}  \sinh(\pi/2) \over \sqrt{\cos(\delta) (1+\cos(\delta))} }\sum_{m=0}^\infty |a_m|  e^{-\delta m}  < \infty  $$
for $\delta \in \left[0, \  \pi/2 \right)$. Hence we proved (2.5). In the same manner the inversion formula (2.6) can be established for the $Im$-transform (1.10). 

\end{proof}

On the other hand, considering discrete transformations with the incomplete modified Bessel functions (1.17), (1.18) 

$$f(x)= \sum_{n=0}^\infty a_n \ {\rm Re} J(x,  1/2+in, \pi),\quad x >0,\eqno(2.9)$$

$$f(x)= \sum_{n=0}^\infty a_n \ {\rm Im} J(x,  1/2+in, \pi),\quad x >0,\eqno(2.10)$$
we have 

{\bf Theorem 3.} {\it Let the sequence $\{a_n\}_{n\ge 0}$ satisfy the $l_1$-condition

$$  \sum_{n=0}^\infty |a_n| < \infty.\eqno(2.11)$$
Then transformations $(2.9), (2.10)$ have the reciprocal inversions, correspondingly,

$$a_n= {4\over \pi^2} \cosh(\pi n) \int_0^\infty  {\rm Re}  K_{1/2+in}(x)  f(x) dx,\eqno(2.12)$$

$$a_n= {4\over \pi^2} \cosh(\pi n) \int_0^\infty  {\rm Im}  K_{1/2+in}(x)  f(x) dx,\eqno(2.13)$$
where integrals $(2.12), (2.13)$ converge absolutely.}

\begin{proof} In fact, since via (1.8), (2.11)

$$\int_0^\infty \left| \left\{\begin{aligned} {\rm Re} \\ {\rm Im} \end{aligned}\right\}   K_{1/2+in}(x)  f(x) \right| dx \le  e^{-\delta n} \int_0^\infty  e^{-x}  K_{1/2}(x \cos(\delta)) dx   $$

$$\int_0^\pi  \cosh( u/2) du \sum_{m=1}^\infty | a_m | = {\pi \sqrt{2}  \sinh(\pi/2)  e^{-\delta n} \over \sqrt{\cos(\delta) (1+\cos(\delta))} }  \sum_{m=0}^\infty | a_m | < \infty $$
when  $\delta \in \left[0, \  \pi/2 \right)$, we have, recalling (1.14) for $\alpha=1/2$,

$${4\over \pi^2} \cosh(\pi n) \int_0^\infty  {\rm Re}  K_{1/2+in}(x)  f(x) dx =  {2\over \pi}    \int_0^\pi \cos(nu)  \sum_{m=0}^\infty a_m  \cos(mu) du = a_n.$$
This gives (2.12).  A companion formula (2.13) can be done immediately, substituting $f$ by  (2.10) into the right-hand side of (2.13), changing the order of integration and summation and using (1.15) for $\alpha =1/2$.   

\end{proof}

Concerning discrete transformations (1.11), (1.12), their inversion formulas will be proved by the same method,  being employed in  [5].  Indeed, it gives

{\bf Theorem 4}. {\it Let $f$ be a complex-valued function on $\mathbb{R}_+$ which is represented by the integral 

$$f(x)= \int_{-\pi}^\pi e^{-x\cosh(u)}  \varphi(u)  du,\ x >0,\eqno(2.14)$$
where $ \varphi(u) = \psi(u)\cosh(u/2)$ and $\psi$ is a  $2\pi$-periodic function, satisfying the Lipschitz condition on $[-\pi, \pi]$, i.e.

$$\left| \psi(u) - \psi(v)\right| \le C |u-v|, \quad  \forall \  u, v \in  [-\pi, \pi]\eqno(2.15)$$
and  $C >0$ is an absolute constant.   Then for all $x >0$ the following inversion formula holds for transformation $(1.11)$

$$f(x)= {a_0\over 2}\   {\rm Re} J(x, 1/2, \pi) +   \sum_{n=1}^\infty \cosh(\pi n)  {\rm Re} J(x,  1/2+in, \pi) a_n.\eqno(2.16)$$
If, in turn,   $ \varphi(u) = \psi(u)\sinh(u/2)$ and $\psi$ is a  $2\pi$-periodic function, satisfying the Lipschitz condition $(2.15)$ on $[-\pi, \pi]$, then 
for all $x >0$ the inversion formula  for transformation $(1.12)$

$$f(x)=   \sum_{n=1}^\infty \cosh(\pi n)  {\rm Im} J(x,  1/2+in, \pi) a_n\eqno(2.17)$$
takes place. }

\begin{proof}  Employing  (1.11), (1.14) for $\alpha=1/2$ and (2.14), we write a  partial sum $S_N(x)$  of the series  (2.16) in the form

$$S_N(x) =   \sum_{n=1}^N \cosh(\pi n)  {\rm Re} J(x,  1/2+in, \pi) a_n$$

$$ =  {2\over \pi} \sum_{n=1}^N  {\rm Re} J(x,  1/2+in, \pi)  \int_{-\pi}^\pi {\varphi(u) \over \cosh(u/2) } \cos(nu) du,\eqno(2.18)$$
where the interchange of the order of integration is, evidently,  permitted by Fubini's theorem.  Moreover, recalling representation (1.17) of the incomplete modified Bessel function and calculating the finite sum with the use of the known formula, we derive

$$S_N(x) =  {1\over \pi} \sum_{n=1}^N  \int_{-\pi}^\pi e^{-x\cosh(t)} \cosh(t/2) \cos(n t) dt  \int_{-\pi}^\pi {\varphi(u) \over \cosh(u/2) } \cos(nu) du $$

$$= {1\over 4\pi}  \int_{-\pi}^\pi e^{-x\cosh(t)} \cosh(t/2)    \int_{-\pi}^\pi {\varphi(u) + \varphi(-u) \over \cosh(u/2) }   {\sin \left((2N+1) (u-t)/2 \right)\over \sin( (u-t) /2)} du dt$$

$$-  {1\over 2\pi}  \int_{-\pi}^\pi e^{-x\cosh(t)} \cosh(t/2) dt    \int_{-\pi}^\pi {\varphi(u)  \over \cosh(u/2) } \ du $$

$$= {1\over 4\pi}  \int_{-\pi}^\pi e^{-x\cosh(t)} \cosh(t/2)    \int_{-\pi}^\pi  \left[ \psi(u) + \psi(-u)\right] {\sin \left((2N+1) (u-t)/2 \right)\over \sin( (u-t) /2)} du dt$$

$$-   {a_0\over 2}\   {\rm Re} J(x, 1/2, \pi).\eqno(2.19)$$
But since $\psi$ is $2\pi$-periodic, we treat  the latter integral with respect to $u$ as follows 

$$  \int_{-\pi}^{\pi}  \left[ \psi(u)+ \psi(-u) \right]  \  {\sin \left((2N+1) (u-t)/2 \right)\over \sin( (u-t) /2)}  du $$

$$=  \int_{ t-\pi}^{t+ \pi}  \left[ \psi(u)+\psi(-u) \right]  \  {\sin \left((2N+1) (u-t)/2 \right)\over \sin( (u-t) /2)}  du $$

$$=  \int_{ -\pi}^{\pi}  \left[ \psi(u+t)+ \psi(-u-t) \right]  \  {\sin \left((2N+1) u/2 \right)\over \sin( u /2)}  du.\eqno(2.20) $$
Moreover,

$$ {1\over 2\pi} \int_{ -\pi}^{\pi}  \left[ \psi(u+t)+ \psi(-u-t) \right]  \  {\sin \left((2N+1) u/2 \right)\over \sin( u /2)}  du - \left[ \psi(t)+ \psi(-t) \right] $$

$$=  {1\over 2\pi} \int_{ -\pi}^{\pi}  \left[ \psi(u+t)- \psi(t) - \psi (-t) + \psi(-u-t) \right]  \  {\sin \left((2N+1) u/2 \right)\over \sin( u /2)}  du.$$
When  $u+t > \pi$ or  $u+t < -\pi$ then we interpret  the value  $\psi(u+t)- \psi(t)$ by  formulas

$$\psi(u+t)+\psi(t) = \psi(u+t-2\pi)+ \psi(t - 2\pi),$$ 

$$\psi(u+t)+ \psi(t) = \psi(u+t+ 2\pi)+ \psi(t +2\pi),$$ 
respectively.  Analogously, the value  $\psi(-u-t)+ \psi(-t)$  can be treated.   Then   due to the Lipschitz condition (2.15) we have the uniform estimate
for any $t \in [-\pi,\pi]$

$${\left|  \psi(u+t)- \psi(t) - \psi (-t) + \psi(-u-t) \right| \over | \sin( u /2) |}  \le 2C \left| {u\over \sin( u /2)} \right|.\eqno(2.21)$$
Therefore,  owing to the Riemann-Lebesgue lemma

$$\lim_{N\to \infty } {1\over 2\pi} \int_{ -\pi}^{\pi}  \left[ \psi(u+t)+ \psi(-u-t)  - \psi(t) - \psi (-t) \right]  \  {\sin \left((2N+1) u/2 \right)\over \sin( u /2)}  du =  0,\eqno(2.22)$$
i.e.

$$\lim_{N\to \infty } {1\over 2\pi} \int_{ -\pi}^{\pi}  \left[ \psi(u+t)+ \psi(-u-t)  \right]  \  {\sin \left((2N+1) u/2 \right)\over \sin( u /2)}  du =  \psi(t) + \psi (-t)\eqno(2.23)$$
for all $ t\in [-\pi,\pi].$    Moreover, returning to (2.19),  and taking into account (2.14), (2.23), it yields

$$\lim_{N\to \infty }  {1\over 4\pi}  \int_{-\pi}^\pi e^{-x\cosh(t)} \cosh(t/2)    \int_{-\pi}^\pi \left[ \psi(u) + \psi(-u)\right] {\sin \left((2N+1) (u-t)/2 \right)\over \sin( (u-t) /2)} du dt$$

$$=   {1\over 2}  \int_{-\pi}^\pi e^{-x\cosh(t)} \cosh(t/2)  \left[  \psi(t) + \psi (-t) \right] dt = \int_{-\pi}^\pi e^{-x\cosh(t)} \varphi(t) dt = f(x).\eqno(2.24)$$
Indeed, we have via (2.20), (2.22)

$$ {1\over 4\pi}  \int_{-\pi}^\pi e^{-x\cosh(t)} \cosh(t/2)    \int_{-\pi}^\pi  \left[ \psi(u) + \psi(-u)\right] {\sin \left((2N+1) (u-t)/2 \right)\over \sin( (u-t) /2)} du dt$$

$$-   {1\over 2}  \int_{-\pi}^\pi e^{-x\cosh(t)} \cosh(t/2)  \left[  \psi(t) + \psi (-t) \right] dt $$

$$=  {1\over 4\pi}  \int_{-\pi}^\pi e^{-x\cosh(t)} \cosh(t/2)    \int_{-\pi}^\pi  \left[ \psi(u+t) + \psi(-u-t)\right.$$

$$\left. -  \psi(t) - \psi (-t) \right] {\sin \left((2N+1) (u-t)/2 \right)\over \sin( (u-t) /2)} du dt \to 0,\quad N \to \infty$$
owing to the dominated convergence theorem since (see (2.21))

$$ \int_{-\pi}^\pi e^{-x\cosh(t)} \cosh(t/2)  \int_{ -\pi}^{\pi} \Bigg| \left[ \psi(u+t)- \psi(-u-t)  - \psi(t) + \psi (-t) \right]\Bigg.$$

$$ \left.  {\sin \left((2N+1) u/2 \right)\over \sin( u /2)}  \right| du dt \le  2 C \int_{-\pi}^\pi e^{-x\cosh(t)}  \cosh(t/2) dt   \int_{ -\pi}^{\pi}   \left| {u\over \sin( u /2)} \right| du < \infty.$$
Thus, passing to the limit in (2.18),  taking into account (2.19), (2.24), we obtain

$$\lim_{N\to \infty } S_N(x) =  f(x)- {a_0\over 2}\   {\rm Re} J(x, 1/2, \pi).$$
This proves (2.16). In order to prove (2.17), we employ (1.12), (1.15), (1.18) to write the partial sum of the series (2.17) in the form

$$ S_N(x) =   \sum_{n=1}^N \cosh(\pi n)  {\rm Im} J(x,  1/2+in, \pi) a_n$$

$$ =  {2\over \pi} \sum_{n=1}^N  {\rm Im} J(x,  1/2+in, \pi)  \int_{-\pi}^\pi {\varphi(u) \over \sinh(u/2) } \sin(nu) du$$

$$=  {1\over 4\pi}  \int_{-\pi}^\pi e^{-x\cosh(t)} \sinh(t/2)    \int_{-\pi}^\pi  \left[ \psi(u) - \psi(-u)\right] {\sin \left((2N+1) (u-t)/2 \right)\over \sin( (u-t) /2)} du dt.$$
Then the same ideas and estimates are applied to get the equality $\lim_{N\to \infty } S_N(x) =  f(x)$,  which yields  (2.17), completing the proof of Theorem 4. 

\end{proof} 

Finally, let us consider analogs of the discrete Lebedev-Skalskaya transforms (1.11), (1.12) with the incomplete modified  Bessel functions (1.17), (1.18)

$$a_n =  {4\over \pi^2} \int_0^\infty  {\rm Re} J(x,  1/2+in, \pi) f(x)dx,\quad n \in \mathbb{N}_0,\eqno(2.25)$$

$$a_n =   {4\over \pi^2} \int_0^\infty  {\rm Im} J(x,  1/2+in, \pi) f(x) dx,\quad n \in \mathbb{N}_0.\eqno(2.26)$$

We have

{\bf Theorem 5}. {\it Let $f$ be represented by series $(1.9)$

$$f(x)= \sum_{m=0}^\infty b_m  {\rm Re} K_{\alpha+im}(x),\eqno(2.27)$$
where the sequence $\{b_m\}_{m\ge 0}$ satisfies condition $(2.1)$. Then transformation $(2.25)$ has the reciprocal inversion formula in the form

$$f(x)= \sum_{n=0}^\infty a_n \cosh(\pi n)  {\rm Re} K_{\alpha+in}(x), \ x > 0,\eqno(2.28)$$
where the series converges absolutely.  Besides, if $f$ is given by series $(1.10)$

$$f(x)= \sum_{m=1}^\infty c_m  {\rm Im} K_{\alpha+im}(x)\eqno(2.29)$$
under the same condition for $\{c_m\}_{m\ge 1}$, then transformation $(2.26)$ is inverted by the formula

$$f(x)= \sum_{n=1}^\infty a_n \cosh(\pi n)  {\rm Im} K_{\alpha+in}(x), \ x > 0\eqno(2.30)$$
with the absolute convergence of the series $(2.30)$.}

\begin{proof} The proof is straightforward via substitution  of (2.25) into the right-hand side of (2.28), using (2.27) and definition of the incomplete modified Bessel function (1.17). Then changing the order of summation and integration due to  the estimate (2.4), we calculate the inner integral by (1.14) and use the orthogonality of trigonometric functions. Thus we arrive at (2.28).  In the same manner one proves (2.30), invoking (2.26), (2.29) and (1.15).

\end{proof}

\bigskip
\centerline{{\bf Acknowledgments}}
\bigskip

\noindent The work was partially supported by CMUP, which is financed by national funds through FCT (Portugal)  under the project with reference UIDB/00144/2020.

\bigskip
\centerline{{\bf References}}
\bigskip
\baselineskip=12pt
\medskip
\begin{enumerate}

\item[{\bf 1.}\ ]  D.S. Jones, Incomplete Bessel functions. I,  {\it Proc. Edinb. Math. Soc.} {\bf 50} (2007), N 1,  173-183.

\item[{\bf 2.}\ ]  N.N. Lebedev, I.P. Skalskaya,  Some integral transforms related to the Kontorovich-Lebedev transform,   {\it Problems in Math. Physics, Leningrad}, 1976, 68-79 (in Russian). 

\item[{\bf 3.}\ ] A.P. Prudnikov, Yu.A. Brychkov and O.I. Marichev, {\it Integrals and Series}. Vol. I: {\it Elementary
Functions}, Vol. II: {\it Special Functions}, Gordon and Breach, New York and London, 1986, Vol. III : {\it More special functions},  Gordon and Breach, New York and London,  1990.

\item[{\bf 4.}\ ] S. Yakubovich, {\it Index Transforms}, World Scientific Publishing Company, Singapore, New Jersey, London and
Hong Kong, 1996.

\item[{\bf 5.}\ ]  S. Yakubovich, Discrete Kontorovich-Lebedev transforms, {\it The Ramanujan Journal} (to appear).

\end{enumerate}

\vspace{5mm}

\noindent S.Yakubovich\\
Department of  Mathematics,\\
Faculty of Sciences,\\
University of Porto,\\
Campo Alegre st., 687\\
4169-007 Porto\\
Portugal\\
E-Mail: syakubov@fc.up.pt\\

\end{document}